# ATTRIBUTING A PROBABILITY TO THE SHAPE OF A PROBABILITY DENSITY

By Peter Hall and Hong Ooi

*Australian National University*

We discuss properties of two methods for ascribing probabilities to the shape of a probability distribution. One is based on the idea of counting the number of modes of a bootstrap version of a standard kernel density estimator. We argue that the simplest form of that method suffers from the same difficulties that inhibit level accuracy of Silverman's bandwidth-based test for modality: the conditional distribution of the bootstrap form of a density estimator is not a good approximation to the actual distribution of the estimator. This difficulty is less pronounced if the density estimator is oversmoothed, but the problem of selecting the extent of oversmoothing is inherently difficult. It is shown that the optimal bandwidth, in the sense of producing optimally high sensitivity, depends on the widths of putative bumps in the unknown density and is exactly as difficult to determine as those bumps are to detect. We also develop a second approach to ascribing a probability to shape, using Müller and Sawitzki's notion of excess mass. In contrast to the context just discussed, it is shown that the bootstrap distribution of empirical excess mass is a relatively good approximation to its true distribution. This leads to empirical approximations to the likelihoods of different levels of "modal sharpness," or "delineation," of modes of a density. The technique is illustrated numerically.

**1. Introduction.** Assigning a probability, or a measure of likelihood, to a quantity determined by an infinite number of unknown parameters is an intrinsically difficult problem. This is particularly the case when definition of the function requires a certain level of smoothing, for example in the case of a probability density. It has recently been proposed [Efron and Tibshirani (1998)] that relative likelihoods of the numbers of modes of a density might be calculated by, in effect, counting the numbers of modes of bootstrap versions of a kernel density estimator. This can be viewed as a development of Silverman's (1981) bootstrap method for testing for the number of modes of a distribution; Silverman adjusted the bandwidth of the estimator until the mode count agreed with that specified by the null hypothesis.

We argue that such bootstrap likelihoods do not converge in probability and in particular do not converge to the "truth" in the standard frequentist sense, unless the bandwidth is chosen an order of magnitude larger than would be appropriate for standard kernel density estimation. Using subsampling methods does not overcome this difficulty; if anything, those techniques make matters a little worse. The difficulties are related to the known level inconsistency of Silverman's (1981) test for the number of modes. Indeed, both problems are rooted in the fact that the bootstrap distribution of a kernel density estimator is not a good approximation to the unconditional distribution of the estimator, if the bandwidth is of its usual pointwise optimal size.

If the bandwidth is allowed to take larger than usual values, then these problems recede. However, the difficulty then arises of determining how large the bandwidth should be. We show that this problem is essentially insoluble; the size of the bandwidth depends on the widths of small potential modes, the very existence of which one is trying to determine.

There is, however, a second, related class of problems, where we may exploit the fact that (under the assumption of a given number of modes) the "modal sharpness," or extent of delineation of the modes of a density, can be accurately estimated in terms of empirical excess mass. Most important, in contrast to problems related to the likelihood of the number of modes, the distribution of empirical excess mass can be accurately approximated using the standard bootstrap, without requiring choice of a smoothing parameter. In this way a set of graphs of constrained density estimates can be constructed, having excess masses that correspond to quantiles of the estimated distribution of excess mass for a given number ($k$, say) of modes, and actually having $k$ modes. A value of $k$ can be determined by testing, or sets of graphs can be constructed for different numbers of modes.

We discuss these two approaches as much because they contrast as because they are similar. The first, density estimator-based technique cannot be interpreted in frequentist terms, and indeed Donoho's (1988) results essentially imply that there is not a meaningful way of empirically assessing the likelihood that a density has $m$ modes. The second method sidesteps these difficulties by eschewing the problem of computing a likelihood for modality, and instead focuses on measuring the "pointiness" of the density's peaks and troughs. Since this approach has a conventional interpretation in frequentist terms, then it is necessarily different from the first, but the two are plainly connected; the number of modes of a probability density is closely related to its excess mass, not least through the fact that empirical approximations to the latter are used to test hypotheses about the former.

Another similarity is that both methods claim to attribute a probability to the shape of a density. For the first method the probability is the likelihood that the density has a given number of modes, while for the second it is the coverage probability of a confidence interval for excess mass.

Section 2.1 will discuss issues of bandwidth choice in the first class of problems, where the likelihood of modality is approximated by mode counting. The second problem class, where shape is described in terms of excess mass, will be discussed in Section 2.2. Both accounts rely critically on theoretical properties, which will be described in Section 3. Efron and Tibshirani (1998) have already carefully worked through numerical examples in the first problem class, and so in the numerical work in Section 2.2 we shall confine attention to the second class.

## 2. Methodology and general properties.

2.1. *Counting modes of a kernel estimator.* Efron and Tibshirani (1998) developed an engaging and particularly original approach to solving problems that are more general than that considered in the present paper. Efron and Tibshirani's method allows them to ascribe a "probability" to the event that a density $f$ is bimodal, by in effect counting the number of modes in the bootstrap form of an appropriately constructed kernel estimator, for example,

$$(2.1) \qquad \hat{f}(x) = \frac{1}{nh} \sum_{i=1}^{n} K\left(\frac{x - X_i}{h}\right),$$

where $K$ is a known probability density, $h$ denotes the bandwidth and $\mathcal{X} = \{X_1, \ldots, X_n\}$ is a random sample drawn from the distribution with density $f$. The arguments Efron and Tibshirani use are, of necessity because the range of problems they treat is so broad, heuristic rather than rigorously mathematical. We shall argue that, in the context of modality of densities, their definition of the amount of probability attributable to different density shapes is not interpretable as a probability in the usual frequentist sense.

Related issues were addressed by Donoho (1988), who demonstrated on essentially topological grounds that the probability that a density has at least $k$ modes is definable, whereas the chance that the density has exactly $k$ modes is not. We should mention too that Efron and Tibshirani's method is somewhat more sophisticated than the counting approach we shall discuss below, for example through being founded on Gaussian-based transformations of the bootstrap distributions of numbers of counts. However, since the large sample distributions of those counts are not approximately Normally distributed (see, e.g., Theorem 3.1), it is not difficult to show that our conclusions apply to the more complex method.

The method that Efron and Tibshirani (1998) suggest using for bandwidth selection, that is, ten-fold cross-validation, produces (as it is designed to) a bandwidth of an order that is asymptotically optimal for pointwise accuracy of the estimator. In particular, the bandwidth is of size $n^{-1/5}$, where $n$ denotes sample size. This will prove important in the first part of our discussion, although later we shall consider larger bandwidths.

We shall show in Theorem 3.1 that for such a bandwidth, the bootstrap distribution of the number of modes of the bootstrapped density estimator converges in distribution but not in probability; the latter, not the former, is the usual sense in which, for practical reasons, one wishes a bootstrap quantity to converge. Moreover, the in-distribution limit does not accurately reflect the number of modes of the sampled distribution. In particular, even if the sampled distribution is strictly unimodal, the weak limit of the bootstrap likelihood is nondegenerately supported on the set of all strictly positive integers.

Naturally one seeks a way of overcoming these difficulties. The method of subsampling, or the "$m$ out of $n$ bootstrap" as it is sometimes called, has a good reputation for remedying convergence problems in a wide range of applications of the bootstrap. See, for example, Bickel and Ren (1996), Lee (1999) and Politis, Romano and Wolf (1999). In Theorem 3.2 we shall show, however, that in the context of estimating the number of modes of $f$, subsampling actually tends to impair performance of the bootstrap when the bandwidth is of size $n^{-1/5}$. It results in the likelihood being approximated by an indicator function, and so the bootstrap estimate of the probability that the sampled density has $k$ modes is well approximated by a random variable that takes only the values 0 and 1. This indicator variable does not converge in probability. It does converge in distribution, but not to the deterministic indicator of the number of modes of the true density.

The landscape changes markedly when a larger order of bandwidth is employed, however. If modes and local minima of the density are "clearly defined," in the sense that the curvature of the density does not vanish at those turning points and the density has no shoulders, and if the bandwidth converges to 0 at a strictly slower rate than $n^{-1/5}$, then the probability that the density estimator has the same number of modes as the true density converges to 1 as $n \to \infty$. This result, and those discussed earlier, are valid provided we avoid spurious small modes in the tails that arise from data sparseness. This problem is commonly addressed as part of kernel-based inference for the number of modes; see, for example, Fisher, Mammen and Marron (1994) and Hall and York (2001). [A density $f$ has a "shoulder" at a point $x$ if both $f'(x)$ and $f''(x)$ vanish and $x$ is not a turning point. To remove the latter possibility it is usual to assume $f''(x) \neq 0$ if $f'(x) = 0$.]

In reality, however, the issues are more complex than this simple asymptotic account suggests. The most important problems involving determination of the number of modes are arguably those where the modes are not "clearly defined" in the context discussed immediately above. Examples include problems where it is difficult to distinguish between a small mode and a shoulder. To some extent these instances too can be satisfactorily addressed by simply counting the number of modes of a kernel density estimator, as proposed by Efron and Tibshirani (1998), although now the choice of bandwidth becomes a more critical issue. Theorem 3.3 will show that the bandwidth should now be at least an order of magnitude larger than $n^{-1/7}$; otherwise, spurious additional modes will be introduced in the region of a shoulder, if the density should have a shoulder rather than a small mode. Therefore, the bandwidth for the density estimator that will enable a bump to be detected must be strictly narrower than the bump for which we are looking.

Moreover, the bandwidth should not be too large, or we shall smooth the bump into the shoulder and miss it altogether. For example, suppose a small bump is constructed above the shoulder, of width $h_1$ and with its height chosen so that the density estimator continues to have three bounded derivatives. Take $h_1 = h_1(n)$ to converge to 0 as $n \to \infty$, so that the problem becomes more complex as more information becomes available. In order to correctly distinguish the bump as a mode, by counting the number of modes of a kernel density estimator, the bandwidth for the latter must converge to zero at a rate that is strictly faster than $h_1$. These results, which are made concise in Theorem 3.4, also hold if we count the number of modes of the bootstrap form of the density estimator.

2.2. *Excess mass as a descriptor of density shape.* The notion of excess mass was introduced by Müller and Sawitzki (1991), and has been discussed extensively; see, for example, Polonik (1995, 1998), Gezeck, Fischer and Timmer (1997), Cheng and Hall (1998), Chaudhuri and Marron (1999, 2000), Polonik and Yao (2000) and Fisher and Marron (2001). It is closely related to Hartigan and Hartigan's (1985) notion of a "dip" in a distribution, and in fact Hartigan and Hartigan's dip test for unimodality is equivalent, in one dimension, to the excess mass test. Either approach can be thought of as being based on the "taut string" method for constructing an empirical distribution that is constrained to be unimodal. That technique has a range of applications to other problems, including monotone and convex approximation [e.g., Leurgans (1982)], nonparametric regression more generally [e.g., Mammen and van de Geer (1997) and Davies and Kovac (2001)] and data exploration [Davies (1995)].

Excess mass of order $m \geq 1$, and the corresponding excess mass difference, are defined, respectively, by

$$E_m(\lambda) = \sup_{L_1,\ldots,L_m} \sum_{i=1}^{m} \{F(L_i) - \lambda \|L_i\|\},$$

(2.2)

$$\Delta_m = \sup_{\lambda > 0} \{E_m(\lambda) - E_{m-1}(\lambda)\},$$

in which the first supremum is taken over all sequences $L_1, \ldots, L_m$ of disjoint intervals, and $\|L\|$ denotes the length of $L$. The empirical form, $\widehat{\Delta}_m$, of $\Delta_m$ is obtained by replacing $E_m(\lambda)$ and $E_{m-1}(\lambda)$ at (2.2) by $\widehat{E}_m(\lambda)$ and $\widehat{E}_{m-1}(\lambda)$, respectively, where $\widehat{E}_m(\lambda)$ is defined as at (2.2) but with $F$ replaced by the empirical distribution function $\widehat{F}$ based on the dataset $\mathcal{X}$. Properties of $\widehat{\Delta}_m$ directly reflect those of $\Delta_m$, not least through the fact that $\widehat{\Delta}_m$ is consistent for $\Delta_m$ as $n \to \infty$.

To appreciate the connection between $\Delta_m$ and the shape of the density $f$, observe that when $m = 2$ and $f$ is bimodal, $\Delta_m$ equals the least amount of mass that needs to be removed from one of the modes, and placed into the trough between them, in order to render $f$ unimodal. In particular, $\Delta_2 = 0$ if and only if $f$ is unimodal. For a general $f$ and for $m \geq 2$, $\Delta_m = 0$ if $f$ has no more than $m - 1$ modes, although the converse is not generally true for $m > 3$. For instance, $\Delta_3 = 0$ for a strictly trimodal density if and only if the height of either of the outer modes does not exceed the height of the local minimum between the other two modes. One can reasonably argue that in this case the lowest mode is insubstantial relative to the other two, and that "$\Delta_3 = 0$ if and only if the density $f$ has no more than two relatively substantial modes." Analogous interpretations are valid for $m \geq 4$.

The fact that the excess mass statistic does not exactly relate to the number of modes (not least because "insubstantial modes" do not directly influence the statistic) means that our approach to ascribing a probability to density shape is quite different from the mode-focused method suggested by Efron and Tibshirani (1998). Our approach is clearly influenced by modality, but is far from being driven by it. It measures the shape of a distribution using information about mode "strength," and in some ways pays scant attention to the number of modes. It is partially linked to Efron and Tibshirani's (1998) approach through work of Chaudhuri and Marron (1999, 2000), which emphasizes mode counts but nevertheless assesses the strengths of putative modes.

Focusing on the case $\Delta_m = 0$ addresses only one example of the ways in which $\Delta_m$ reflects the shape of $f$. More generally, the fact that $E_m(\lambda)$ represents the maximum deviation of $f$ from a composition of $m$ uniform distributions of height $\lambda$ implies that as $\Delta_m$ increases at least some of the modes of $f$ become more pronounced. To gain insight into this property, consider the case where $f$ is the density of a mixture of $p \geq m$ Normal $N(\mu_i, \sigma_i^2)$ populations, with distinct fixed means $\mu_1, \ldots, \mu_p$ and respective nonvanishing, fixed mixing proportions $\pi_1, \ldots, \pi_p$. The supremum of $\Delta_m$, over all such densities, equals the sum of the $p - m + 1$ smallest values of $\pi_i$, and is attained by letting the corresponding $p - m + 1$ values of $\sigma_i$ decrease to zero. In particular, the modes corresponding to these distributions in the mixture become infinitely sharp spikes. Furthermore, if the mixture density has $p$ modes then the maximum value of $\Delta_m$ is attained only by letting the $p - m + 1$ values of $\sigma_i$ (corresponding to the $p - m + 1$ smallest $\pi_i$'s) converge to zero.

The simplest case in the Normal mixture example is that where $m = p = 2$ and $0 < \pi_1 < \pi_2 < 1$. There, $\Delta_2 > 0$ if and only if $\sigma_1$ and $\sigma_2$ are chosen so that the Normal mixture is bimodal, and $\sigma_1 \to 0$ as $\Delta_2$ increases to its maximum value, $\pi_1$. (The limit $\Delta_2 \to \pi_1$ can be attained with $\sigma_2$ fixed and $\sigma_1 \to 0$, and also with $\sigma_1$, $\sigma_2 \to 0$ together.)

Properties such as those discussed in the three previous paragraphs argue that for fixed $m$, relatively large values of $\Delta_m$ are associated with densities $f$ that have more than $m - 1$ "substantive" modes, and with all but $m - 1$ of the modes being relatively sharp.

2.3. *Imposing constraints on excess mass.* If we were to construct $\hat{f}$ in such a way that $\Delta_m(\hat{f}) = \widehat{\Delta}_m$ then we would be allowing the estimator to reflect the actual empirical level of "modal sharpness." Note particularly that calculation of $\widehat{\Delta}_m$ does not involve any smoothing, whereas $\hat{f}$ does require a smoothing parameter. Conceptually, computing $\hat{f}$ subject to $\Delta_m(\hat{f}) = \widehat{\Delta}_m$ is similar to constructing a density estimator subject to one or more moments of the distribution with density $\hat{f}$ being equal to the corresponding empirical moments for the dataset $\mathcal{X}$. Advantages of the latter procedure have been discussed by Jones ([1991](#)) and Hall and Presnell ([1999](#)), for example. In practice, however, the task is significantly more difficult when the constraint is in terms of excess mass, not least because $\Delta_m(\hat{f})$ is a highly nonlinear function of $\hat{f}$.

We can be more bold than to ask simply that $\Delta_m(\hat{f}) = \widehat{\Delta}_m$. The distribution of $\widehat{\Delta}_m$ may be approximated using bootstrap methods (see Theorem 3.5), and estimates of the quantiles of the distribution may be computed. In this way we may construct versions of $\hat{f}$ under the constraint that its excess mass equals any given quantile, thereby computing density estimates that reflect the sharpness of the true density in a median sense or in the sense of any given probability for excess mass.

This procedure can be implemented using data-sharpening methods [Choi and Hall ([1999](#)) and Braun and Hall ([2004](#))], to impose constraints on estimator shape. The method produces a new estimator $\hat{f}_{\mathcal{Y}}$, computed as was $\hat{f} = \hat{f}_{\mathcal{X}}$ but from a sharpened dataset $\mathcal{Y}$, with excess mass $\hat{\Delta}$, say. [We would usually choose $\hat{\Delta}$ to be an estimator of a quantile of the distribution of $\Delta_m(\hat{f}_{\mathcal{Y}})$.] The method starts with a density estimator, $\bar{f}$, which could be either $\hat{f}_{\mathcal{X}}$ or $\hat{f}_{\mathcal{Z}}$, the latter being another version of $\hat{f}_{\mathcal{X}}$, this time constructed for another sharpened sample $\mathcal{Z}$.

In the latter case, $\mathcal{Z}$ might be deliberately constructed so that $\hat{f}_{\mathcal{Z}}$ has a different shape from $\hat{f}_{\mathcal{X}}$. Discussion of the principle of data sharpening, and of reasons why an intermediate dataset, $\mathcal{Z}$, might be generated from $\mathcal{X}$ prior to using data sharpening to impose a constraint on excess mass, is given in Appendix A. An algorithm for data sharpening is presented in Appendix B.

The number of modes of $\hat{f}_{\mathcal{Y}}$ will be determined partly by the number of modes of $\bar{f}$, and partly by the numerical value chosen for $\hat{\Delta}$. For example, if $\bar{f}$ is unimodal, implying that $\Delta_2(\bar{f}) = 0$, but we take $m = 2$ and $\hat{\Delta} > 0$, then $\hat{f}_{\mathcal{Y}}$ will have a second mode, generally becoming more pronounced as $\hat{\Delta}$ increases. If $\bar{f}$ is trimodal then constraining $\Delta_3(\hat{f}_{\mathcal{Y}})$ to equal $\hat{\Delta} < \Delta_3(\bar{f})$, and steadily reducing $\hat{\Delta}$ to zero, may reduce the number of modes to two or may simply reduce the height of one of the two outer modes to the height of the local minimum between the other two modes. The outcome here depends on $\bar{f}$. If $\bar{f}$ is trimodal, and if one of the modes is "insubstantial" (in the sense of Section 2.2), then constraining $\Delta_3(\hat{f}_{\mathcal{Y}})$ to equal $\hat{\Delta} > \Delta_3(\bar{f})$ will often make that mode more pronounced; the mode will not be smoothed away. (More generally, all the modes to which we refer above are modes in the usual sense, i.e., local maxima of the density. They can be either "substantial" or "insubstantial" from the viewpoint of excess mass.)

There are potential alternative approaches, although they are difficult to implement in practice. One might consider using a single, fixed bandwidth, and vary it to ensure a given value of empirical excess mass. However, this approach is so strongly influenced by data in the tails of the distribution that it is often impractical. For example, if $f$ is a Normal density, and the desired number of modes equals one, then the bandwidth must diverge to infinity with sample size, at rate at least $(\log n)^{1/2}$, in order to ensure that empirical excess mass difference equals zero (equivalently, that the density estimator is unimodal). An alternative technique would be to use a bandwidth that varies with location, but that approach too is strongly influenced by outlying data and is difficult to use to estimate densities with a given number of modes when that number exceeds one. Moreover, even under the constraint of a single mode it is difficult to select a variable bandwidth that produces a given value of excess mass.

2.4. *Real-data illustration of constraints on excess mass.* We illustrate application of our data-sharpening method to the chondrite dataset of Good and Gaskins (1972, 1980). There is evidence [e.g., Leonard (1978) and Silverman (1981)] that these data come from a distribution with at least two modes, and likely no more than two modes. Good and Gaskins (1972), Simonoff (1983) and Minnotte, Marchette and Wegman (1998) note that the chondrite dataset may contain evidence of three modes, while Müller and Sawitzki (1991) are inconclusive in this regard. Evidence for the third mode is based on just three data points, and so is not strong; see Silverman's contribution to the discussion of Leonard (1978). A kernel density estimate based on the chondrite data, with bandwidth chosen by the Sheather and Jones (1991) plug-in rule, is shown by the dotted lines in the panels of Figure 1 and does in fact have just two modes. If it had three modes, say, we would use data sharpening to reduce one of the modes to a shoulder, so that the final density estimate had just two modes. Then in subsequent steps of our algorithm we would replace the real dataset by its sharpened form.

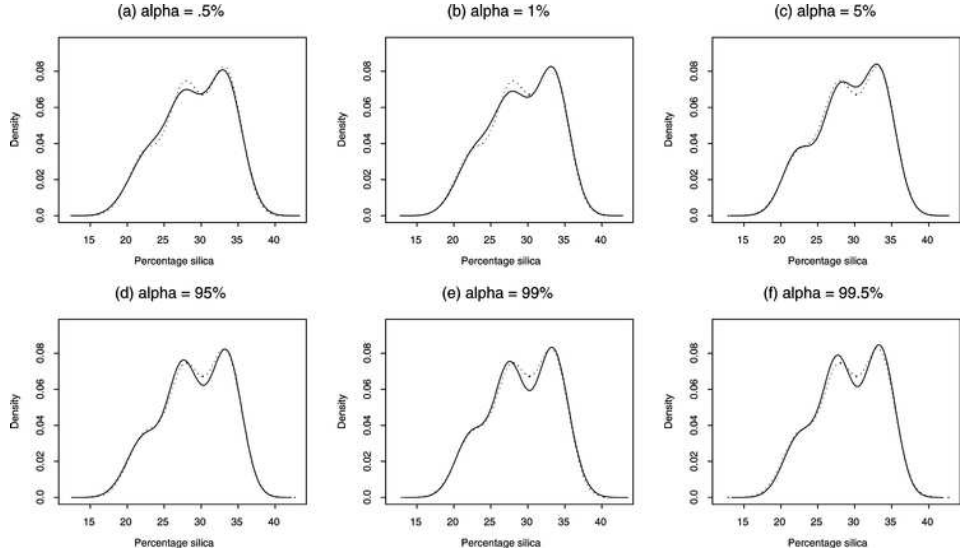

Fig. 1. *Density estimates calculated from sharpened versions of the chondrite data, where the extent of sharpening is such as to ensure the excess mass of the density estimate equals the $\alpha$-level quantile of the bootstrap distribution of the excess mass statistic. In each panel the dotted curve depicts the conventional kernel estimator, using the same Sheather–Jones bandwidth as the sharpened versions shown by the unbroken curve. The values of $\alpha$ are 0.005, 0.01, 0.05, 0.95, 0.99 and 0.995, and correspond to the curves shown in panels (a)–(f), respectively.*

Fixing the bandwidth, and using bootstrap simulation, we estimated quantiles of the distribution of $\widehat{\Delta}_2$ for levels $\alpha = 0.005, 0.01, 0.05, 0.95, 0.99$ and $0.995$. The corresponding density estimates are depicted in Figure 1. They were computed using the algorithm given in Appendix B. Estimates for low values of $\alpha$ are relatively close to being unimodal, while those for $\alpha$ close to 1 have pronounced modes and antimodes.

We also applied our technique to the geyser dataset of Weisberg (1985) and Scott (1992), which consists of 107 eruption durations for the Old Faithful geyser. Tests for multimodality based on the excess mass statistic, calibrated in a variety of ways, argue strongly that the sampled distribution has at least two modes; see Müller and Sawitzki (1991). There is no evidence of more than two modes, and in particular a kernel density estimate constructed using the Sheather and Jones (1991) plug-in bandwidth shows two pronounced modes and not even a suggestion of a shoulder. Construction of the quantile curve estimates gives results similar to those in Figure 1.

## 3. Theoretical properties of shape probabilities.

3.1. *Shape probabilities for fixed densities.* Let $\hat{f}$ be as at (2.1), and denote by $\hat{f}^*$ the standard bootstrap form of $\hat{f}$, computed from a resample $\mathcal{X}^*$ derived by sampling randomly with replacement from $\mathcal{X}$. Write $N(n)$ and $N^*(n)$ for the numbers of modes (i.e., local maxima) of $\hat{f}$ and $\hat{f}^*$, respectively. Among other results we shall show that if $n^{1/5}h \to C_0 > 0$ as $n \to \infty$ then $N(n)$ and $N^*(n)$ converge in distribution. The limit is degenerate if and only if $C_0 = \infty$; in this case it is concentrated at the atom 1.

Next we describe the limiting distributions of $N(n)$ and $N^*(n)$ in the "standard" case, where $n^{1/5}h \to C_0 \in (0, \infty)$. Let $W$ and $W^*$ denote independent standard Brownian bridges, let $x_0$ be the mode of $f$, and assuming $f(x_0) > 0$ and $f''(x_0) < 0$, put

$$\xi(y) = f(x_0)^{1/2} \int K''(u) W(y+u) \, du,$$

$$\xi^*(y) = f(x_0)^{1/2} \int K''(u) W^*(y+u) \, du,$$

(3.1)

$$\eta(y) = C_0^{-3/2} \xi(y) + C_0 y f''(x_0),$$

$$\eta^*(y) = C_0^{-3/2} \{\xi(y) + \xi^*(y)\} + C_0 y f''(x_0),$$

each stochastic process being defined for $-\infty < y < \infty$. Note that when $K$ is the Gaussian kernel, each process is infinitely differentiable with probability 1. Let $N$ and $N^*$ denote the numbers of downcrossings of 0 by $\eta$ and $\eta^*$, respectively. Both random variables are well defined and finite with probability 1 and take only strictly positive integer values. We shall note in Theorem 3.1 that, under regularity conditions, the limiting distribution of $N(n)$ is the distribution of $N$, and the limit of the bootstrap distribution of $N^*(n)$ may be expressed as the distribution of $N^*$ conditional on $W$.

Assume $f$ has two continuous derivatives on its support, which we take to equal $\mathcal{S} = [a, b]$ where $-\infty < a < b < \infty$, and that $f(a) = f(b) = 0$, $f'(a+) > 0$ and $f'(b-) < 0$. Call this condition ($C_{f1}$). Suppose too that in the interior of $\mathcal{S}$ the equation $f'(x) = 0$ has a unique solution $x_0 \in (a, b)$, and that $f''(x_0) < 0$; call this $(C_{f2})$. Assume of the bandwidth that for some $\delta > 0$, $h = h(n) = O(n^{-\delta})$ as $n \to \infty$, and that $n^{1/5}h$ is bounded away from 0; call this condition $(C_h)$. For simplicity, and since the Gaussian kernel is by far the most commonly used in density estimation problems associated with shape, we shall suppose throughout that $K(u) = (2\pi)^{-1/2}\exp(-u^2/2)$. However, since monotonicity of the number of modes of $\hat{f}$ as a function of $h$ is not a concern in our work, the majority of our results hold for sufficiently smooth, unimodal, compactly supported kernels such as the triweight. In such cases the inequalities $\hat{f}' > 0$ and $\hat{f}' < 0$ in the third probability at (3.3) should be replaced by nonsharp inequalities.

Let $\hat{x}_0$ denote the point at which $\hat{f}$ achieves its largest local maximum. Then $\hat{x}_0$ is well defined with probability 1.

THEOREM 3.1 (Bootstrap approximation to distribution of $N$). *Assume* $(C_{f1})$, $(C_{f2})$ *and* $(C_h)$, *and that* $K$ *is the Gaussian kernel.*

(a) *If in addition to* $(C_h)$ *we have* $n^{1/5}h \to C_0 \in (0, \infty)$, *then*

$$\sup_{k \geq 0} |P\{N(n) = k\} - P(N = k)| \to 0$$

*as* $n \to \infty$. *While this result continues to hold if* $N(n)$ *and* $N$ *are replaced by* $N^*(n)$ *and* $N^*$, *respectively, the bootstrap distribution of* $N^*(n)$ *does not converge in the usual sense. Indeed, there exists a construction of* $(W, W^*)$ *that depends on* $\mathcal{X}$ *and is such that*

$$(3.2) \qquad \sup_{k \geq 0} |P\{N^*(n) = k|\mathcal{X}\} - P(N^* = k|W)| \to 0$$

*in probability as* $n \to \infty$.

(b) *If, on the other hand* $n^{1/5}h \to \infty$, *then both* $P\{N(n) = 1\}$ *and* $P\{N^*(n) = 1\}$ *converge to 1, and so with probability converging to 1 both* $\hat{f}$ *and* $\hat{f}^*$ *are unimodal. Furthermore, if* $n^{(1/5)-\delta}h \to \infty$ *for some* $\delta > 0$, *then each of the probabilities*

$$(3.3) \qquad P\{N(n) = 1\}, \qquad P\{N^*(n) = 1\} \quad and$$

$$P\{\hat{f}' > 0 \ on \ (-\infty, \hat{x}_0) \ and \ \hat{f}' < 0 \ on \ (\hat{x}_0, \infty)\}$$

*equals* $1 - O(n^{-\lambda})$ *for all* $\lambda > 0$.

The first portions of parts (a) and (b) of Theorem 3.1, relating only to the nonbootstrap case, are given by Mammen (1995). See also Mammen, Marron and Fisher (1992) and Konakov and Mammen (1998).

It will follow from our proof of (3.2) that the particular construction of $W$, given the data, does not converge, and in particular that $P(N^* =$ $k|W)$ does not converge in probability as $n \to \infty$. Therefore, when $C_0 < \infty$ the distribution of $N^*(n)$, conditional on the data, does not converge in probability as $n \to \infty$. Furthermore, part (a) of the theorem implies that while the unconditional distribution of $N^*(n)$ does converge, it does not converge to the limiting distribution of $N(n)$. Theorem 3.1 has several more general or more detailed forms, which are given in a longer version of this paper obtainable from the authors.

Next we show that subsampling fails to remove the inconsistency problems suffered by the bootstrap in the present setting. In fact the bootstrap distribution of $N^*(n)$, in the case of subsampling, is well approximated by a crude indicator function of $N(n)$. Nevertheless, subsampling does not, to first order, impair consistency when the bandwidth is of larger order than $n^{-1/5}$. These properties are stated formally in Theorem 3.2. By way of notation, we redefine $\mathcal{X}^* = \mathcal{X}^*(m)$ to be a resample of size $m < n$ drawn by sampling randomly, with replacement, from $\mathcal{X}$, and construct $\hat{f}^*$ and $N_n^*$ for this version of $\mathcal{X}^*$.

THEOREM 3.2 (Subsample bootstrap approximation to distribution of $N$). *Assume* $(C_{f1})$, $(C_{f2})$ *and* $(C_h)$, *that* $K$ *is the Gaussian kernel, and that the resample size* $m = m(n)$ *satisfies* $m \to \infty$ *and* $m/n \to 0$ *as* $n \to \infty$.

(a) *If in addition to* $(C_h)$ *we have* $n^{1/5}h \to C_0 \in (0, \infty)$, *then*

$$\sup_{k \geq 0} |P\{N^*(n) = k | \mathcal{X}\} - I\{N(n) = k\}| \to 0$$

*in probability as* $n \to \infty$.

(b) *If on the other hand* $n^{1/5}h \to \infty$, *then both* $P\{N(n) = 1\}$ *and* $P\{N^*(n) = 1\}$ *converge to* 1.

3.2. *Shape probabilities for densities with small modes.* Theorem 3.1 has analogues in the case of densities with one or more shoulders, that is, points at which $f'$ and $f''$ both vanish. In this case the critical size of the bandwidth is $n^{-1/7}$, rather than $n^{-1/5}$, provided $f'''$ does not also vanish at the shoulder.

In particular, if $h$ is of smaller order than $n^{-1/7}$ then the probability that the number of modes of $\hat{f}$ exceeds any fixed integer converges to 1 as $n \to \infty$, and if $n^{1/7}h \to C_0 > 0$ then the distribution of the number of modes has a proper limit, degenerate at the atom $\nu$ if $f$ has just $\nu$ modes and $C_0 = \infty$. For the latter result it is sufficient to assume $(C_{f1})$, along with the condition $(C_{f4})$ that in the interior of the support of $f$ the equation $f'(x) = 0$ has just $\omega$, say, solutions, at just $2\nu - 1$ of which $f'' \neq 0$, with $f'''$ having three continuous, nonvanishing derivatives in the neighborhoods of the other $\omega - 2\nu + 1$ zeros of $f'$. Constraint $(C_{f4})$ implies that $f$ has $\nu$ local maxima, $\nu - 1$ local minima and $\omega - 2\nu + 1$ shoulders.

Shoulders may be regarded as embryonic modes, and for this reason densities with shoulders are of particular interest since they lie on boundaries separating classes of densities with different shapes, expressed through their "modalities." See, for example, Cheng and Hall (1999). In both theoretical and numerical studies the performances of methods for assigning probabilities to density shapes may be assessed in terms of their success in distinguishing between densities that have shoulders and those which have small modes in places that would otherwise be shoulders. With this in mind we shall expand the class of densities satisfying $(C_{f1})$ and $(C_{f3})$ by allowing the first and second derivatives, but not the first, second and third, to vanish simultaneously. We shall discuss the performance, uniformly over such densities, of empirical methods for assigning probabilities to the numbers of modes and show that techniques based on counting the number of modes of a kernel density estimator can have optimal performance, in a minimax sense, if bandwidth is chosen larger than $n^{-1/7}$.

To simplify discussion we shall base our lower bound on perturbations of a density $f$ with just one mode and one shoulder. Specifically, $f$ will satisfy $(C_{f1})$ and the following condition, which we call $(C_{f5})$: In the interior of $\mathcal{S}$ the equation $f'(x) = 0$ has just two solutions, $x_0, x_1 \in (a, b)$, with $x_0$ denoting the mode of $f$ and satisfying $f''(x_0) < 0$, and $x_1$ representing a shoulder and such that $f$ has three continuous derivatives in a neighborhood of $f$, $f''(x_1) = 0$ and $f'''(x_1) \neq 0$. Given any empirical procedure $\mathcal{N}$ for counting the number of modes of a density, we would want $\mathcal{N}$ to equal 1, with high probability, when applied to a dataset drawn from a distribution whose density satisfies $(C_{f1})$ and $(C_{f5})$.

Now perturb $f$ by adding a small bump at the shoulder, as follows. Let $\psi$ denote a symmetric, compactly supported probability density with three continuous derivatives on the real line, a unique mode at the origin satisfying $\psi''(0) < 0$, no other point $x$ in the interior of the support of $\psi$ such that $\psi'(x) = 0$ and such that the equation $\frac{1}{2}|f'''(x_1)|y^2 = |\psi'(y)|$ has a unique solution on $(0, \infty)$ which also satisfies $f'''(x_1)y + \psi''(y) \neq 0$. Call this condition $(C_\psi)$; as we shall show in the proof of Theorem 3.3, the last part of $(C_\psi)$ ensures that the added bump produces a single additional mode. Let $h_1 = h_1(n) \to 0$, and let the perturbed density be

$$f_n(x) = \frac{f(x) + h_1^3 \psi\{(x - x_1)/h_1\}}{1 + h_1^4}.$$

In this formula, the factor $h_1^3$ ensures that like $f$, $f_n$ has three bounded derivatives in a neighborhood of $x_1$. The denominator $1 + h_1^4$ guarantees that $f_n$ integrates to 1. We would want $\mathcal{N}$ to equal 2, with high probability, when applied to data from the distribution with density $f_n$.

The density $f_n$ has just two modes, at $y_0 = x_0 + o(h_1)$ and $y_1 = x_1 + o(h_1)$, respectively. Of course, a local minimum occurs between them [at a point with formula $x_1 + O(h_1)$], but no other turning points and no shoulders exist. Choosing $h_1$ larger or smaller makes the small bump near $x_1$ more or less pronounced, respectively.

Our next theorem shows that in order for it to be possible to correctly distinguish two modes in the density $f_n$, based on a sample of size $n$, the rate at which $h_1$ converges to 0 must be strictly slower than $n^{-1/7}$.

THEOREM 3.3 (Necessity of using large bandwidth when counting modes). *Assume $f$ satisfies* (C$_{f1}$) *and* (C$_{f5}$), *and that $\psi$ satisfies* (C$_\psi$). *Let $\mathcal{N}$ denote any empirical procedure for counting the number of modes of a density, and use it to estimate the number of modes of $f$ and of $f_n$, based on samples of size $n$ from these respective distributions. If $\mathcal{N}$ is asymptotically correct in each case, that is, if both $P_f(\mathcal{N} = 1) \to 1$ and $P_{f_n}(\mathcal{N} = 2) \to 1$, then $n^{1/7}h_1 \to \infty$ as $n \to \infty$.*

It is likewise possible to show that if $n^{1/7}h_1 \to \infty$ then, provided the bandwidth $h$ in the kernel density estimator $\hat{f}$ converges to 0 at a rate that lies strictly between $n^{-1/7}$ and the rate at which $h_1$ decreases, the naive rule $\mathcal{N}$ that simply counts the number of modes of $\hat{f}$ is asymptotically correct. In this sense it achieves the level of precision that is shown by Theorem 3.3 to be optimal. See Mammen (1995) for discussion and details. Konakov and Mammen (1998) treat the multivariate version of this result.

3.3. *Probability distribution of excess mass.* In Section 2.2 we defined the excess mass, $\Delta_m$, of $f$, and discussed potential applications of approximations to the distribution of the empirical form, $\hat{\Delta}_m$, of this quantity. Here we describe the limiting distribution of $\hat{\Delta}_m$ in the case $m = 2$.

Given intervals $L_1, L_2$ as at (2.2), let $L_{0j}$ denote the version of $L_j$ that produces the second supremum there in the case $m = 2$. Assume $f = F'$ is bimodal, let $\lambda_0$ denote the value of $\lambda$ that maximises $E_2(\lambda) - E_1(\lambda)$, and write $L_{01} = (x_1, x_2)$ and $L_{02} = (x_3, x_4)$, where without loss of generality, $x_1 < \cdots < x_4$. In this notation,

$$
\begin{aligned}
(3.4) \quad E_2(\lambda_0) = &\{F(x_2) - F(x_1) - \lambda_0(x_2 - x_1)\} \\
&+ \{F(x_4) - F(x_3) - \lambda_0(x_4 - x_3)\}.
\end{aligned}
$$

Note that $f(x_i) = \lambda_0$ for $1 \le i \le 4$.

We shall suppose too that one mode contains strictly less excess mass than the other, in the sense that the mode from which mass is removed, and placed into the trough between the modes when the nearest unimodal density is constructed, is uniquely defined. We shall call this mode the "smallest mode." Without loss of generality the smallest mode is the second of the two modes, lying between $x_3$ and $x_4$. See Figure 5 of Müller and Sawitzki (1991) for an illustration of this case.

Next we define the limiting distribution of $\widehat{\Delta}_2$; it is a mixture of correlated Normals. Let $N = (N_2, N_3, N_4)$ denote a trivariate Normally distributed vector with zero mean and covariances given by $\operatorname{var}(N_i, N_j) = F(x_i)\{1 - F(x_j)\}$ for $i \leq j$; redefine $\xi_1$ and $\xi_2$ to be standard Brownian motions, stochastically independent of $N$; and define $I$ to equal 1 if

$$\frac{\sup\{\xi_2(u) - u^2\}}{\sup\{\xi_4(u) - u^2\}} < \left|\frac{f'(x_2)}{f'(x_4)}\right|^{1/3}$$

and to equal 2 otherwise. Finally, put $Z = N_{2I} - N_3$.

We assume of $f$ that it has a continuous derivative, ultimately monotone in each tail, that the constraints $f'(x) = 0$ and $f(x) \neq 0$ are jointly satisfied at just three points, $x^{(1)} < x^{(2)} < x^{(3)}$, in the neighborhood of each of which $f''$ exists and is continuous, and $f''(x^{(1)}) < 0$, $f''(x^{(2)}) > 0$ and $f''(x^{(3)}) < 0$, and the points $x_1, \ldots, x_4$ at (3.4) are such that each $f'(x_i) \neq 0$. Call this condition $(\mathrm{C}_{f6})$. Let $\widehat{\Delta}_2^*$ denote the version of $\widehat{\Delta}_2$ computed not from $\mathcal{X}$ but from $\mathcal{X}^*$, the latter obtained by sampling randomly, with replacement, from $\mathcal{X}$.

THEOREM 3.4 (Consistency of bootstrap estimate of excess mass distribution). *Assume $f$ satisfies $(\mathrm{C}_{f6})$. Then the distribution of $n^{1/2}(\widehat{\Delta}_2 - \Delta_2)$ converges, as $n \to \infty$, to the distribution of $Z$. Furthermore, the conditional distribution of $n^{1/2}(\widehat{\Delta}_2^* - \widehat{\Delta}_2)$, given $\mathcal{X}$, converges in probability to the distribution of $Z$.*

It follows from Theorem 3.4, and symmetry of the distribution of $Z$, that both the standard percentile bootstrap methods consistently estimate quantiles of the distribution of $\widehat{\Delta}_2$. Therefore, the percentile bootstrap produces confidence intervals for the excess mass $\Delta_2$ that have asymptotically correct coverage accuracy. For instance, if $0 < \alpha < 1$ and $\hat{t}_\alpha$ is defined to be the infimum of values $t$ such that $P(\widehat{\Delta}_2^* \leq t | \mathcal{X}) \geq \alpha$, then it follows from the theorem that $P(\Delta_2 \leq \hat{t}_\alpha) \to \alpha$ as $n \to \infty$. The percentile bootstrap technique was used in Section 2.2 to construct confidence regions for $\Delta_2$ and hence to compute estimates of $f$ whose shapes (in terms of their "modal sharpness" or "delineation," as expressed through excess mass) correspond to particular quantiles.

Theorem 3.4 is readily extended to show that, under regularity conditions analogous to $(\mathrm{C}_{f6})$, and for general $m \geq 2$, the limiting distribution of $n^{1/2}(\Delta_m - \Delta_m)$ is consistently approximated by the conditional distribution $n^{1/2}(\widehat{\Delta}_m^* - \widehat{\Delta}_m)$. Of course, such a result fails if, when computing the bootstrap approximation, we mistakenly constrain the initial estimator $\hat{f}_{\mathcal{Z}}$ to have too few modes. (See Section 2.3 for discussion of $\hat{f}_{\mathcal{Z}}$.) In particular, if we constrain $\hat{f}_{\mathcal{Z}}$ to have $M < m_0$ modes, where $m_0$ denotes the true number of modes of $f$, then $\hat{f}_{\mathcal{Z}}$ converges not to $f$ but to an $M$-mode density that is nearest to $f$ in a sense that can be defined in terms of the distance measure, $d$, used for the data-sharpening algorithm. On this occasion, this basic inconsistency renders invalid any bootstrap approximations that start from $\hat{f}_{\mathcal{Z}}$.

**4. Proofs of Theorems 3.1 and 3.2.** The nonbootstrap parts of Theorem 3.1 are given by Mammen (1995), but since parts of his argument are needed for the bootstrap case and for Theorem 3.2, they are reproduced in outline form here.

4.1. *Monotonicity of $\hat{f}$ outside $(x_0 - \varepsilon, x_0 + \varepsilon)$.* Define $\widehat{D}_1(x) = \hat{f}'(x) - E\{\hat{f}'(x)\}$ and $\ell = \log n$. The argument used to derive Lemma 6 of Mammen, Marron and Fisher (1992) [see also Silverman (1983)] may be employed to prove that for each $\lambda > 0$ there exists $B = B(\lambda) > 0$ such that

$$P\left\{ \sup_{a \le x \le b} |\widehat{D}_1(x)| > B(\ell/nh^3)^{1/2} \right\} = O(n^{-\lambda}).$$

Note too that $E\{\hat{f}'(x)\} = f'(x) + o(1)$ uniformly in $x \in [a + \delta, b - \delta]$ for each $\delta > 0$, while $E\{\hat{f}'(x)\} \ge \frac{1}{2} f'(x) + o(1)$ uniformly in $x \in [a, x_0]$ and $E\{\hat{f}'(x)\} \le \frac{1}{2} f'(x) + o(1)$ uniformly in $x \in [x_0, b]$. It follows from these properties that for each $\varepsilon \in (0, \min(x_0 - a, b - x_0))$, and all $\lambda > 0$,

$$P\{\hat{f}' > 0 \text{ on } [a, x_0 - \varepsilon] \text{ and } \hat{f}' < 0 \text{ on } [x_0 + \varepsilon, b]\} = O(n^{-\lambda}).$$

The definition of $\hat{f}$ implies directly that with probability 1, $\hat{f}' > 0$ on $(-\infty, a)$ and $\hat{f}' < 0$ on $(b, \infty)$. Hence, for each $\varepsilon \in (0, \min(x_0 - a, b - x_0))$, and all $\lambda > 0$,

$$(4.1) \quad \begin{aligned} P\{\hat{f}' > 0 \text{ on } (-\infty, x_0 - \varepsilon] \text{ and } \hat{f}' < 0 \text{ on} &= [x_0 + \varepsilon, \infty)\} \\ &= 1 - O(n^{-\lambda}). \end{aligned}$$

Similarly, (4.1) holds if $\hat{f}'$ is replaced by $(\hat{f}^*)'$.

4.2. *Approximation to $\hat{f}'$ and $(\hat{f}^*)'$ on $(x_0 - \varepsilon, x_0 + \varepsilon)$.* Define $\widehat{D}_1^*(x) = (\hat{f}^*)'(x) - \hat{f}'(x)$,

$$\xi_1(x) = \int K''(u) W_0\{F(x + hu)\} \, du,$$

$$\xi_1^*(x) = \int K''(u) W_0^*\{\widehat{F}(x + hu)\} \, du,$$

where $W_0$, and $W_0^*$ conditional on both $\mathcal{X}$ and $W_0$, are standard Brownian bridges, and $\hat{F}$ denotes the conventional empirical distribution function of the sample $\mathcal{X}$ from which $\hat{f}$ was computed. It may be proved, using the embedding of Komlós, Major and Tusnády (1976), that $W_0$ and $W_0^*$ may be constructed such that

$$(4.2) \qquad \begin{aligned} \hat{D}_1(x) &= n^{-1/2}h^{-2}\xi_1(x) + R_1(x), \\ \hat{D}_1^*(x) &= n^{-1/2}h^{-2}\xi_1^*(x) + R_1^*(x), \end{aligned}$$

where for each $\delta, \lambda > 0$ and $\varepsilon \in (0, \min(x_0 - a, b - x_0))$,

$$(4.3) \qquad \begin{aligned} P\left\{\sup_{|x-x_0|\leq\varepsilon} |R_1(x)| > n^{-1+\delta}h^{-2}\right\} &= O(n^{-\lambda}), \\ P\left\{\sup_{|x-x_0|\leq\varepsilon} |R_1^*(x)| > n^{-1+\delta}h^{-2}\right\} &= O(n^{-\lambda}). \end{aligned}$$

4.3. *Monotonicity of $\hat{f}'$ and $(\hat{f}^*)'$ outside* $(x_0 - Ch, x_0 + Ch)$. Note that $E\{\hat{f}'(x)\} = f'(x) + o(h)$ uniformly in $x \in (x_0 - \varepsilon, x_0 + \varepsilon)$, for sufficiently small $\varepsilon > 0$, and that

$$(4.4) \qquad \sup_{|y|\leq C} |E\{\hat{f}'(x_0 + hy)\} - hyf''(x_0)| = o(h)$$

for any $C > 0$. It may be deduced from these results, the fact that $h$ is not less than a constant multiple of $n^{-1/5}$, and properties of a Brownian bridge, that for each $C_1, \delta > 0$ there exists $C > 0$ such that for all sufficiently small $\varepsilon > 0$, and all sufficiently large $n$,

$$(4.5) \qquad \begin{aligned} &P\{n^{-1/2}h^{-2}\xi_1(x) + E\{\hat{f}'(x)\} \text{ is greater than } C_1 h \\ &\quad \text{for } -\varepsilon \leq x - x_0 \leq -Ch, \\ &\quad \text{and is less than } -C_1 h \text{ for } Ch \leq x - x_0 \leq \varepsilon\} \geq 1 - \delta. \end{aligned}$$

Similarly we may prove that (4.5) holds if we replace $\xi_1$ by $\xi_1 + \xi_1^*$. If $n^{(1/5)-\delta}h \to \infty$ for some $\delta > 0$ then both results may be strengthened by replacing "$\geq 1 - \delta$" on the right-hand side of (4.5) by "$= 1 - O(n^{-\lambda})$ for all $\lambda > 0$."

Combining (4.1)–(4.3) and (4.5) we deduce that for each $\delta > 0$ there exists $C > 0$ such that for all sufficiently large $n$,

$$(4.6) \qquad \begin{aligned} &P\{\hat{f}'(x) \text{ is strictly positive for } x \leq x_0 - Ch \\ &\quad \text{and strictly negative for } x \geq x_0 + Ch\} \geq 1 - \delta, \end{aligned}$$

and that (4.6) continues to hold if $\hat{f}'$ is replaced by $(\hat{f}^*)'$. Moreover, both results continue to hold with "$\geq 1 - \delta$" on the right-hand side of (4.6) replaced by "$= 1 - O(n^{-\lambda})$ for all $\lambda > 0$," provided $n^{(1/5)-\delta}h \to \infty$ for some $\delta > 0$.

4.4. *Approximation to $\hat{f}'$ and $(\hat{f}^*)'$ on $(x_0 - Ch, x_0 + Ch)$.* Define $\widehat{D}(y) = \widehat{D}_1(x_0 + hy)$ and $\widehat{D}^*(y) = \widehat{D}_1^*(x_0 + hy)$. Noting the continuity properties of a Brownian bridge and that for each $\delta, \lambda > 0$ the probability that $\sup |\widehat{F} - F|$ exceeds $n^{\delta - (1/2)}$ equals $O(n^{-\lambda})$, we may deduce that, defining

$$\xi_2^*(y) = \int K''(u) W_0^* \{ F(x_0 + hy + hu) \} \, du,$$

it is true that for each $C, \delta, \lambda > 0$,

$$P \left\{ \sup_{|y| \le C} |\xi_1^*(x_0 + hy) - \xi_2^*(y)| > n^{\delta - (1/4)} \right\} = O(n^{-\lambda}).$$

From this result, (4.2), (4.3) and the fact that $h$ is no smaller than a constant multiple of $n^{-1/5}$, we may deduce that for each $C, \lambda > 0$ and some $\delta > 0$,

$$(4.7) \quad P \left\{ \sup_{|y| \le C} |\widehat{D}^*(y) - n^{-1/2} h^{-2} \xi_2^*(y)| > n^{-(1/2) - \delta} h^{-3/2} \right\} = O(n^{-\lambda}).$$

Note that we may write $W_0(t) = V(t) - tV(t)$, where $V$ is a standard Brownian motion, and that $W_0^*$ may be represented analogously. These properties, and arguments similar to those in the previous paragraph, allow us to show that if we define $\xi$ and $\xi^*$ as at (3.1), for appropriate choices of $W$ and $W^*$, then for some $\delta > 0$ we have for each $C, \lambda > 0$,

$$P \left[ \sup_{|y| \le C} \{ |\xi_1(y) - h^{1/2} \xi(y)| + |\xi_2^*(y) - h^{1/2} \xi^*(y)| \} > n^{-\delta} h^{1/2} \right] = O(n^{-\lambda}).$$

From this result, (4.2), (4.3) and (4.7) we may deduce that

$$(4.8) \quad \begin{aligned} \widehat{D}(y) &= (nh^3)^{-1/2} \{ \xi(y) + R(y) \}, \\ \widehat{D}^*(y) &= (nh^3)^{-1/2} \{ \xi^*(y) + R^*(y) \}, \end{aligned}$$

where for some $\delta > 0$ and each $C, \lambda > 0$,

$$(4.9) \quad P \left[ \sup_{|y| \le C} \{ |R(y)| + |R^*(y)| \} \ge n^{-\delta} \right] = O(n^{-\lambda}).$$

In the notation of (4.8),

$$(4.10) \quad \begin{aligned} \hat{f}'(x_0 + hy) &= (nh^3)^{-1/2} \{ \xi(y) + R(y) \} + E\{ \hat{f}'(x_0 + hy) \}, \\ (\hat{f}^*)'(x_0 + hy) &= (nh^3)^{-1/2} \{ \xi(y) + \xi^*(y) + R(y) + R^*(y) \} \\ &\quad + E\{ \hat{f}'(x_0 + hy) \}. \end{aligned}$$

4.5. *Proof of Theorem* 3.1(a). Define $\eta$ and $\eta^*$, in terms of $\xi$ and $\xi^*$, as at (3.1). It follows from (4.4), (4.9) and (4.10) that

$$\zeta(y) \equiv n^{1/5}\hat{f}'(x_0 + hy) = \eta(y) + o_p(1),$$

$$\zeta^*(y) \equiv n^{1/5}(\hat{f}^*)'(x_0 + hy) = \eta^*(y) + o_p(1),$$

both results holding uniformly in $|y| \leq C$. The processes $\eta$, $\eta^*$, $\zeta$ and $\zeta^*$ are all differentiable, each derivative equals $O_p(1)$ uniformly on $[-C, C]$, and $\zeta' = \eta' + o_p(1)$ and $(\zeta^*)' = (\eta^*)' + o_p(1)$ uniformly on $[-C, C]$, for each $C > 0$. Furthermore, if the downcrossings of 0 by $\eta$ on $[-C, C]$ occur at points $Z_1, \ldots, Z_M$, where $M \leq N$, then (a) $P(M = N) \to 1$ as $C \to \infty$, (b) with probability 1 no $Z_i$ equals $C$ or $-C$, and (c) for each $\varepsilon > 0$,

$$\lim_{\varepsilon \to 0} P\{|\eta'(Z_i)| > \varepsilon \text{ for each } i\} = 1.$$

Together these properties imply that for each $C > 0$,

(4.11)
$$P\{\text{the number of downcrossings of 0 by } \hat{f}'$$
$$\text{on } (x_0 - Ch, x_0 + Ch) \text{ equals}$$
$$\text{the number of downcrossings of 0 by } \eta$$
$$\text{on } (x_0 - Ch, x_0 + Ch)\} \to 1$$

as $n \to \infty$. Similarly, (4.11) holds if $(\hat{f}', \eta)$ is replaced by $((\hat{f}^*)', \eta^*)$. Theorem 3.1(a) follows from (4.6), (4.11) and their bootstrap forms.

4.6. *Proof of Theorem* 3.1(b). Minor modifications of the previous arguments show that when $n^{1/5}h \to \infty$, $P\{N(n) = 1\}$ and $P\{N^*(n) = 1\}$ both converge to 1. Next we prove that when $n^{(1/5)-\delta}h \to \infty$ for some $\delta > 0$, we have for each $\lambda > 0$,

(4.12)
$$P\{N(n) = 1\} = 1 - O(n^{-\lambda}).$$

Similar arguments may be used to obtain the same identity for the other two probabilities at (3.3).

Observe from (4.6) and the comments which immediately follow it that (4.12) will follow if we show that for each $C > 0$,

$$P\{\hat{f}' \text{ has at most one zero in } (x_0 - Ch, x_0 + Ch)\} = 1 - O(n^{-\lambda})$$

for each $\lambda > 0$. This result is in turn implied by: for each $\lambda > 0$,

(4.13)
$$P\{\hat{f}'' \text{ has no zeros in } (x_0 - Ch, x_0 + Ch)\} = 1 - O(n^{-\lambda}).$$

We may establish an analogue, for $\hat{f}''$, of the first parts of (4.9) and (4.10),

$$\hat{f}''(x_0 + hy) = (nh^5)^{-1/2}\{\xi'(y) + S(y)\} + E\{\hat{f}''(x_0 + hy)\},$$

where, for some $\varepsilon > 0$ and each $C, \lambda > 0$,

$$(4.14) \qquad P\left\{ \sup_{|y| \leq C} |S(y)| \geq n^{-\varepsilon} \right\} = O(n^{-\lambda}).$$

The condition $n^{(1/5)-\delta}h \to \infty$, which we are currently assuming, implies that $(nh^5)^{-1/2} = O(n^{-\varepsilon})$ for some $\varepsilon > 0$. From this result, (4.14) and properties of a Brownian motion, we may deduce that for each $\varepsilon > 0$ and all $C, \lambda > 0$,

$$P\left\{ (nh^5)^{-1/2} \sup_{|y| \leq C} |\xi'(y) + S(y)| \geq \varepsilon \right\} = O(n^{-\lambda}).$$

It follows from this property, the fact that $f''(0) < 0$ and the expansion $E\{\hat{f}''(x_0 + hy)\} = f''(x_0) + o(1)$ uniformly in $|y| \leq C$, that the probability that $\hat{f}'' < 0$ throughout $(x_0 - Ch, x_0 + Ch)$ equals $1 - O(n^{-\lambda})$ for all $C, \lambda > 0$. This implies (4.13).

4.7. *Outline proof of Theorem* 3.2. Result (4.6), and the properties noted immediately below it, continue to be valid in the present case. And (4.7) holds in the following form, for the same definition of $\xi_2^*$ as before: for each $C, \lambda > 0$ and some $\delta > 0$,

$$P\left\{ \sup_{|y| \leq C} |\hat{D}^*(y) - n^{-1}m^{1/2}h^{-2}\xi_2^*(y)| > n^{-(1/2)-\delta}h^{-3/2} \right\} = O(n^{-\lambda}).$$

Thus, in place of (4.9) and the second part of (4.10) we may write

$$(\hat{f}^*)'(x_0 + hy) = (nh^3)^{-1/2}\{\xi(y) + o_p(1)\} + hyf''(y)$$
$$= \hat{f}'(x_0 + hy) + o_p\{(nh^3)^{-1/2}\},$$

uniformly in $|y| \leq C$. The argument in Section 4.4 may now be used to show that the probability, conditional on $\mathcal{X}$, that the number of downcrossings of 0 by $(\hat{f}^*)'$ equals the number of downcrossings of 0 by $\eta$, converges to 1 as $n \to \infty$. Likewise, the unconditional probability that the number of downcrossings of 0 by $\eta$ equals the number of downcrossings of 0 by $\hat{f}'$ converges to 1. Part (a) of Theorem 3.2 follows from these properties, and part (b) may be derived similarly.

## 5. Proofs of Theorems 3.3 and 3.4.

5.1. *Proof of Theorem* 3.3. First we show that, under condition $(C_\psi)$, the density $f_n$ has just two modes, one local minimum and no shoulders on its support, for all sufficiently small $h_1$; call this property (P). It will follow that $\mathcal{N}$ is asymptotically correct if it concludes (with probability converging to 1 as $n \to \infty$) that $f$ and $f_n$ have just one and two modes, respectively.

In view of the definition of $f_n$, any turning point of $f_n$ on $(a, b)$ that is not identical to $x_0$ must converge to $x_1$ as $n \to \infty$. Assume without loss of generality that $f'''(x_1) > 0$, and note that

$$(1 + h_1^4)f_n'(x_1 + h_1 y) = f'(x_1 + h_1 y) + h_1^2 \psi'(y)$$
$$= h_1^2 \{ \tfrac{1}{2} f'''(x_1) y^2 + \psi'(y) + o(1) \}$$

as $h_1 \to \infty$, uniformly in $|y| \le C$ for any $C > 0$. These formulas, and the assumption [part of $(C_\psi)$] that the equation

$$(5.1) \qquad\qquad \tfrac{1}{2} f'''(x_1) y^2 = |\psi'(y)|$$

has a unique solution $y_0$ in $(0, \infty)$, imply that any turning point $y_1$ of $f_n$ that converges to $x_1$ as $n \to \infty$, and is not identical to $x_1$, must satisfy

$$(5.2) \qquad\qquad y_1 = x_1 + h_1 y_0 + o(h_1),$$

where $y_1$ is the solution of (5.1). Moreover, since $f'''(x_1)y_0 + \psi''(y_0) \ne 0$, again by virtue of $(C_\psi)$, then the equation $f'(x_1 + h_1 y) + h_1^2 \psi'(y) = 0$ can have no more than one solution $y_1$ satisfying (5.2). It follows that $y_1$ must represent a unique local minimum between $x_0$ and $x_1$, and that (P) holds.

We may view $\mathcal{N}$ as a rule for discriminating between $f$ and $f_n$, determining that the $n$-sample from which $\mathcal{N}$ is computed comes from $f_n$ if $\mathcal{N} = 2$ and comes from $f$ otherwise. If $P_f(\mathcal{N} = 1) \to 1$ and $P_{f_n}(\mathcal{N} = 2) \to 1$, then $\mathcal{N}$ provides asymptotically perfect discrimination, and so, by the Neyman–Pearson lemma, the likelihood ratio rule also provides perfect discrimination. It suffices to show that the latter property implies $n^{1/7} h_1 \to \infty$.

We shall argue by contradiction and show that if $n^{1/7} h$ is bounded as $n \to \infty$ through some infinite sequence, $\mathcal{A}$ say, then the likelihood ratio rule does not provide asymptotically perfect discrimination along the sequence. It may be assumed without loss of generality that $nh \to \infty$ as $n \to \infty$ through $\mathcal{A}$, since otherwise a simple subsidiary argument produces a contradiction.

Observe that, in view of the compact support of $\psi$,

$$(1 + h_1^4)\frac{f_n(x)}{f(x)} = 1 + h_1^3 \frac{1}{f(x)} \psi\left( \frac{x - x_1}{h_1} \right),$$

uniformly in $x \in (a, b)$, as $n \to \infty$ through values in $\mathcal{A}$. The log-likelihood ratio is therefore

$$LR \equiv \sum_{i=1}^n \log\{ f_n(X_i)/f(X_i) \}$$

$$= h_1^3 \sum_{i=1}^n \left\{ \frac{1}{f(X_i)} \psi\left( \frac{X_i - x_1}{h_1} \right) - h_1 \right\}$$

$$-\frac{1}{2}h_1^6\sum_{i=1}^n\left\{\frac{1}{f(X_i)}\psi\left(\frac{X_i-x_1}{h_1}\right)\right\}^2+o_p(1)$$

$$=(nh_1^7)^{1/2}\sigma Z-\frac{1}{2}nh_1^7\sigma^2+o_p(1),$$

where $\sigma^2=(\int\psi^2)/f(x_1)$, the random variable $Z$ is asymptotically standard normal, and the remainders are of the stated orders as $n\to\infty$ through $\mathcal{A}$. Therefore $LR=O_p(1)$ as $n\to\infty$ through values in $\mathcal{A}$, and so it is not possible for the likelihood-ratio test to discriminate, with asymptotic probability 1, against $f_n$ for data from $f$ as $n\to\infty$ through $\mathcal{A}$.

5.2. *Outline proof of Theorem* 3.4. We shall derive only the first, unconditional limit theorem; the second, conditional bootstrap result may be proved similarly. At a key point in the latter proof, where the bootstrap form $W_0^*$ of a Brownian bridge is used in the form of a function of the empirical distribution function $\widehat{F}$ (cf. Section 4.2), we may replace $n^{-1/2}W_0^*(\widehat{F})$ by $n^{-1/2}W_0^*(F)$ and incur an error of only $O_p(n^{-3/4}\log n)$. This is of smaller order than the error of size $n^{-2/3}$ that arises if the approximation is subsequently pursued using arguments developed below, in the nonbootstrap case. In this way it can be seen that the "in distribution" limits are identical in the two cases.

Using the embedding of Komlós, Major and Tusnády (1976) we may, for each $n$, construct a standard Brownian bridge $W_0$ such that

$$\widehat{F}(x)=F(x)+n^{-1/2}W_0\{F(x)\}+O_p(n^{-1}\ell),$$

uniformly in $x$, where $\ell=\log n$. Of course, $W_0(t)=B(t)-tB(1)$ for a standard Brownian motion $B$. Put $\eta=n^{-1/3}$, write $y_i=x_i+\eta u_i$ where $\sup_i|u_i|\le C$ for some fixed $C>0$, and define $N_i=W_0\{F(x_i)\}$ and

$$W_i(t)=(\lambda_0\eta)^{-1/2}[B\{F(x_i)+\lambda_0\eta t\}-N_i].$$

Then $W_i$ is a standard Brownian motion, and $F(y_i)=F(x_i)+\lambda_0\eta u_i+O(\eta^2)$. Therefore, using properties of the modulus of continuity of $B$, we deduce that

$$(5.3)\qquad\widehat{F}(y_i)-F(y_i)=n^{-1/2}N_i+(\lambda_0\eta/n)^{1/2}W_i(u_i)+O_p(\eta\ell/n^{1/2})$$

uniformly in $|u_i|\le C$.

Put $\delta_i=\eta u_i$ and define $\Delta$ to denote the operator describing the perturbation arising when $x_i$ is changed to $y_i=x_i+\delta_i$, for $1\le i\le4$, small $|\delta_i|$ and $\lambda=\lambda_0$ held fixed. For example, $\Delta_iF(x_i)=F(x_i+\delta_i)-F(x_i)$. Then, since each $f(x_i)$ equals $\lambda_0$, $\Delta\{F(x_i)-\lambda_0x_i\}=\frac{1}{2}\delta_i^2f'(x_i)+o(\eta^2)$. From this result and (5.3) we deduce that if $\lambda=\lambda_0+\eta^2v$ then

$$\widehat{F}(y_i)-\widehat{F}(y_j)-\lambda(y_i-y_j)$$

$$= F(x_i) - F(x_j) - \lambda_0(x_i - x_j) + \tfrac{1}{2}\eta^2\{f'(x_i)u_i^2 - f'(x_j)u_j^2\}$$
$$+ n^{-1/2}(N_i - N_j) + (\lambda_0\eta/n)^{1/2}\{W_i(u_i) - W_j(u_j)\}$$
$$- \eta^2 v(x_i - x_j) + O_p(\eta\ell n^{-1/2} + \eta^3),$$

uniformly in $|u_i|, |v| \leq C$. Equivalently, if we define the intervals $L = (y_j, y_i)$ and $L_0 = (x_j, x_i)$ then

$$\widehat{F}(L) - \lambda\|L\| - \{F(L_0) - \lambda_0\|L_0\|\}$$
$$= n^{-1/2}(N_i - N_j)$$
$$+ \eta^2[\lambda_0^{1/2}\{W_i(u_i) - W_j(u_j)\} + \tfrac{1}{2}\{f'(x_i)u_i^2 - f'(x_j)u_j^2\} - v(x_i - x_j)]$$
$$+ O_p(\eta\ell n^{-1/2}).$$

Therefore, if $L^{(1)} = (y_1, y_2)$, $L^{(2)} = (y_3, y_4)$, $L_0^{(1)} = (x_1, x_2)$ and $L_0^{(2)} = (x_3, x_4)$ then

$$\sum_{i=1}^{2}\{\widehat{F}(L^{(i)}) - \lambda\|L^{(i)}\|\}$$
$$= E_2(\lambda_0) + n^{-1/2}(N_2 + N_4 - N_1 - N_3)$$
$$+ \eta^2[\lambda_0^{1/2}\{W_2(u_2) + W_4(u_4) - W_1(u_1) - W_3(u_3)\}$$
$$+ \tfrac{1}{2}\{f'(x_2)u_2^2 + f'(x_4)u_4^2 - f'(x_1)u_1^2 - f'(x_3)u_3^2\}$$
$$- v(x_2 + x_4 - x_1 - x_3)] + o_p(\eta^2).$$

Taking the supremum over $u_1, \ldots, u_4$, and noting that $C$ in the bound $|u_i| \leq C$ is an arbitrary although fixed number, we deduce that

$$\begin{aligned}(5.4)\qquad \widehat{E}_2(\lambda) &= E_2(\lambda_0) + n^{-1/2}(N_2 + N_4 - N_1 - N_3)\\ &\quad + \eta^2\lambda_0^{1/2}\sum_{i=1}^{4}\sup_u\{B_i(u) - b_iu^2\}\\ &\quad - \eta^2 v(x_2 + x_4 - x_1 - x_3) + o_p(\eta^2),\end{aligned}$$

where $B_i(u) = (-1)^i W_i(u)$ is a standard Brownian motion process, and $b_i = (-1)^{i+1}f'(x_i) > 0$. Strictly speaking, $\widehat{E}_2(\lambda)$ on the left-hand side and the supremum on the right-hand side are defined with the suprema taken only over $|u_i| \leq C$. However a subsidiary argument shows that (5.4) holds when the suprema are interpreted over the whole real line.

Similarly,

$$\widehat{E}_1(\lambda) = E_1(\lambda_0) + n^{-1/2}(N_{2J} - N_1)$$

$$(5.5) \qquad + \eta^2 \lambda_0^{1/2} \left[ \sup_u \{ B_{2J}(u) - b_{2J} u^2 \} + \sup_u \{ B_1(u) - b_1 u^2 \} \right]$$

$$- \eta^2 v(x_{2J} - x_1) + o_p(\eta^2),$$

where $J = 1$ or $2$ according as

$$\sup_u \{ B_2(u) - b_2 u^2 \} > \sup_u \{ B_4(u) - b_4 u^2 \}$$

is true or false. Subtracting (5.5) from (5.4), taking the supremum over $v$ and writing $I = 3 - J$, we deduce that

$$(5.6) \qquad \widehat{\Delta}_2 = \Delta_2 + n^{-1/2}(N_{2I} - N_3) + O_p(\eta^2).$$

[Much as in the cases of (5.4) and (5.5), a subsidiary argument shows that the suprema over $v$ may be taken over the whole positive real line, not just over $|v| \leq C$.]

We may write $b_i^{1/3} B_i(u) = \xi_i(t)$ where $t = b_i^{2/3} u$ and $\xi_i$ is, like $B_i$, a standard Brownian motion. In this notation, $B_i(u) - b_i u^2 = b_i^{-1/3} \{ \xi_i(t) - t^2 \}$, and so $I = 1$ or $2$ accordingly as

$$\frac{\sup \{ \xi_2(u) - u^2 \}}{\sup \{ \xi_4(u) - u^2 \}} < \left| \frac{f'(x_2)}{f'(x_4)} \right|^{1/3}$$

is true or false, respectively. The variables $N_2$, $N_3$, and $N_4$ have a joint Normal distribution with zero mean and covariances given by $\mathrm{var}(N_i, N_j) = F(x_i)\{1 - F(x_j)\}$ for $i \leq j$. Furthermore, the $N_i$'s are asymptotically independent of the processes $\xi_i$. The theorem follows from these properties and (5.6).

## APPENDIX A

**Description of data sharpening for constraining excess mass.** The method is based on a density estimator, which we shall denote by $\bar{f}$. This can be either a conventional estimator, $\hat{f}$, computed from $\mathcal{X}$, and which we could denote by $\hat{f}_{\mathcal{X}}$ to indicate that fact or an estimator computed after $\mathcal{X}$ has been sharpened to $\mathcal{Z} = \{ Z_1, \ldots, Z_n \}$, say. In this case we denote the estimator by $\hat{f}_{\mathcal{Z}}$. The dataset $\mathcal{Z}$ might be chosen so that $\hat{f}_{\mathcal{Z}}$ has a given number of modes. See the next paragraph for further discussion. Of course, the case $\bar{f} = \hat{f}_{\mathcal{Z}}$ subsumes $\bar{f} = \hat{f}_{\mathcal{X}}$ as a special, degenerate case, so we may take $\bar{f} = \hat{f}_{\mathcal{Z}}$ below.

If one of our aims is to ensure that $\hat{f}_{\mathcal{Z}}$ has just $m$ (say) modes, where $m$ is different from the number of modes of $\hat{f}_{\mathcal{X}}$, then we might proceed as follows. First, choose a bandwidth for the density estimator (usually by employing a standard method applied to the original dataset), and let $d(\cdot, \cdot)$ denote a nonnegative measure of distance on the real line. It need not be a metric, but for ease of interpretation it should be symmetric. For example, $d(x, y) = (x - y)^2$ is a possibility. Put $d(\mathcal{X}, \mathcal{Z}) = \sum_i d(X_i, Z_i)$, and choose $\mathcal{Z}$ to minimise $d(\mathcal{X}, \mathcal{Z})$ subject to the constraint that $\hat{f}_{\mathcal{Z}}$ just has $m$ modes. [See Hall and Kang (2002) for discussion.]

The $m$th mode will in fact be a shoulder, but can be made more pronounced (once the shoulder is achieved) by transferring the constraint to one on excess mass, rather than on the number of modes. Specifically, once the estimator $\hat{f}_{\mathcal{Z}}$ with just $m$ modes is attained, sharpen $\mathcal{Z}$ to $\mathcal{Y}$ by minimizing $d(\mathcal{Z}, \mathcal{Y})$ subject to $\hat{f}_{\mathcal{Y}}$ having an increased value of excess mass; that is, $\Delta_m(\hat{f}_{\mathcal{Y}}) = \hat{\Delta}$, where $\hat{\Delta} > \Delta_m(\hat{f}_{\mathcal{Z}})$ would typically be chosen to be an estimator of a quantile of the distribution of $\Delta_m(f)$.

In each case the constraints may be imposed using methods based on simulated annealing; this approach is elementary in terms of code, although lengthy from the viewpoint of computing time. The algorithm is described in Appendix B.

For example, if we wished to use $m = 3$ in the algorithm discussed above, but the estimator $\hat{f}_{\mathcal{X}}$ had only one mode, the first step would generally be to sharpen $\mathcal{X}$ to $\mathcal{Z}$ so that $\hat{f}_{\mathcal{Z}}$ had three modes. Nevertheless, although the value of $m$ could be greater or less than the actual number of modes of $\hat{f}_{\mathcal{X}}$, usually it would be less than that number, reflecting the fact that standard kernel density estimators (with appropriately chosen bandwidths) tend to have more, not fewer, modes than the true density. Note too that perhaps not all the modes will be substantial, in the sense of excess mass (see Section 2.2).

## APPENDIX B

**Algorithm for data sharpening subject to constraints on excess mass.** Let the "starting" dataset be $\mathcal{Z} = \{Z_1, \ldots, Z_n\}$, and denote by $\hat{\Delta}_m^{(\alpha)}$ our bootstrap estimator of the $\alpha$-level quantile of the excess mass distribution. Let the sharpened dataset be $\mathcal{Y} = \{Y_1, \ldots, Y_n\}$ and define the distance between $\mathcal{Z}$ and $\mathcal{Y}$ to be $D(\mathcal{Z}, \mathcal{Y}) = \sum_i (Z_i - Y_i)^2$. We seek $\mathcal{Y}$ to minimize $D(\mathcal{Z}, \mathcal{Y})$ subject to $\Delta_m(\hat{f}_{\mathcal{Y}}) = \hat{\Delta}_m^{(\alpha)}$. This problem is solved by a standard simulated annealing algorithm, the perturbations of which (within the annealing loop) are generated as follows.

Let $\tilde{f}$ denote the density estimator $\hat{f}_{\mathcal{Z}}$, and write $\tilde{f}_{\max}$ for its maximum value. At the next step of the algorithm we decide whether we wish to make the data less or more "diffuse," based on whether the current excess-mass statistic $\Delta_m(\hat{f}_{\mathcal{Y}})$ is less than or greater than the target value $\hat{\Delta}_m^{(\alpha)}$, respectively. For a given data point $y_i$, where $1 \leq i \leq n$, we generate a move as

(B.1)                     $$y_i \leftarrow y_i + sz_i \exp\{-\tilde{f}(y_i)/\tilde{f}_{\max}\}$$

if we wish to make the data less diffuse or

$$(B.2) \qquad\qquad y_i \leftarrow y_i + sz_i \exp[\{\tilde{f}(y_i) - \tilde{f}_{\max}\}/\tilde{f}_{\max}]$$

if we wish to make them more diffuse. Here $s$ is a constant equal to the range of $\mathcal{Z}$ divided by 1000 (a value which was chosen by trial and error), and $z_i$ is a number drawn randomly from the standard Normal distribution. Using formulas (B.1) and (B.2) to govern the perturbations was found to give better convergence rates than employing a naive perturbation formula.

The perturbation of $\mathcal{Y}$ indicated by (B.1) or (B.2), for $1 \leq i \leq n$, was ignored if it took $\Delta_m(\hat{f}_{\mathcal{Y}})$ further from the target value $\hat{\Delta}_m^{(\alpha)}$. The algorithm was terminated when $\Delta_m(\hat{f}_{\mathcal{Y}})$ got within $s$ of $\hat{\Delta}_m^{(\alpha)}$. We repeated this procedure 100 times and selected as the solution the configuration with the lowest value of $D(\mathcal{Z}, \mathcal{Y}_j)$.

In practice, this algorithm always converged. In numerical experiments, to check whether the limit was significantly affected by early steps taken by the algorithm, we sometimes started it from small perturbations of $\mathcal{Z}$, but nevertheless reached the same limit.

**Acknowledgment.** We are grateful to three reviewers for helpful comments.

CENTRE FOR MATHEMATICS AND ITS APPLICATIONS
AUSTRALIAN NATIONAL UNIVERSITY
CANBERRA, ACT 0200
AUSTRALIA
E-MAIL: Peter.Hall@maths.anu.edu.au